\documentclass[12pt]{article}
\usepackage{latexsym,amsfonts,amsthm,amsmath,amscd,amssymb}
\usepackage[dvips]{graphicx}

\newtheorem{lemma}{Lemma}[section]
\newtheorem{thm}[lemma]{Theorem}
\newtheorem{rem}[lemma]{Remark}
\newtheorem{prop}[lemma]{Proposition}

\newtheorem{conj}[lemma]{Conjecture}

\newtheorem{defn}[lemma]{Definition}

\newcommand\matS{{\mathbb{S}}}

\newcommand\matC{{\mathbb{C}}}
\newcommand\matE{{\mathbb{E}}}
\newcommand\matH{{\mathbb{H}}}
\newcommand\matN{{\mathbb{N}}}

\newcommand\Sigmatil{{\widetilde\Sigma}}

\renewcommand{\hbar}{{\overline{h}}}

\newfont{\Got}{eufm10 scaled 1200}
\newcommand{\permu}{{\hbox{\Got S}}}
\newcommand{\compo}{\,{\scriptstyle\circ}\,}

\newcommand{\mycap} [1] {\caption{\footnotesize{#1}}}

\newcommand{\myto}{\mathop{\longrightarrow}\limits}
\newcommand\Xtil{{\widetilde X}}
\newcommand\ptil{{\widetilde p}}
\newcommand{\chiorb}{\chi^{{\mathrm{orb}}}}
\newcommand{\myprod}{\!\cdot\!}

\newcommand{\dotsto}{\mathop{\dashrightarrow}\limits}

\newcommand{\argdotsto}[2]{\mathop{\dashrightarrow}\limits^{#1}_{\scriptscriptstyle{#2}}}
\newcommand{\argdotstobis}[2]{\mathop{\dashrightarrow\dashrightarrow}\limits^{#1}_{\scriptscriptstyle{#2}}}
\newcommand{\argdotstoter}[2]{\mathop{\dashrightarrow\dashrightarrow\dashrightarrow}\limits^{#1}_{\scriptscriptstyle{#2}}}
\newcommand{\argdotstoqua}[2]{\mathop{\dashrightarrow\dashrightarrow\dashrightarrow\dashrightarrow}\limits^{#1}_{\scriptscriptstyle{#2}}}
\newcommand{\argdotstoqui}[2]{\mathop{\dashrightarrow\dashrightarrow\dashrightarrow\dashrightarrow\dashrightarrow}\limits^{#1}_{\scriptscriptstyle{#2}}}
\newcommand{\argdotstosex}[2]{\mathop{\dashrightarrow\dashrightarrow\dashrightarrow\dashrightarrow\dashrightarrow\dashrightarrow}\limits^{#1}_{\scriptscriptstyle{#2}}}

\begin{document}

\title{Realizability and exceptionality of candidate surface branched covers:
methods and results}

\author{Ekaterina~{\textsc Pervova}\thanks{Supported by the Marie Curie fellowship MIF1-CT-2006-038734}
\and\addtocounter{footnote}{5} Carlo~{\textsc Petronio}}

\maketitle

\begin{abstract}
\noindent
Let $f:\Sigmatil\myto^{d:1}_{\Pi}\Sigma$
denote a branched cover between closed, connected, and orientable surfaces,
where $d$ is the global degree, $\Pi=\{\Pi_1,\ldots,\Pi_n\}$,
$n$ is the number of branching points, and $\Pi_i$ is the partition of
$d$ given by the local degrees over the $i$-th branching point.
If $\ell(\Pi)$ is the total length of $\Pi$, then
the Riemmann-Hurwitz formula asserts that
$\chi(\Sigmatil)-\ell(\Pi)=d\cdot(\chi(\Sigma)-n)$.
A \emph{candidate branched cover} is a symbol $\Sigmatil\argdotsto{d:1}{\Pi}\Sigma$
satisfying the same condition, and it is called \emph{realizable} if
there exists
a corresponding $f:\Sigmatil\myto^{d:1}_{\Pi}\Sigma$.
The problem of determining which candidate covers are realizable is
very old and still partially unsolved.
In this paper we will review five different techniques employed
in recent years to attack it, and we will state the main
results obtained using them. Each technique will be exemplified through
a proof that the candidate cover
$S^2\argdotstoter{4:1}{(2,2),(2,2),(3,1)}S^2$ is \emph{exceptional},
namely non-realizable. We will also state some results for the non-orientable
version of the problem, none of which is due to us.
\end{abstract}

\section{Introduction}
A \emph{branched cover} is a map $f:\Sigmatil\to\Sigma$, where
$\Sigmatil$ and $\Sigma$ are closed connected surfaces and $f$ is
locally modelled on maps of the form $\matC\ni z\mapsto z^k\in\matC$
for some $k\geqslant1$. The integer $k$ is called the \emph{local
degree} at the point of $\Sigmatil$ corresponding to $0$ in the
source $\matC$. If $k>1$ then the point of $\Sigma$ corresponding to
$0$ in the target $\matC$ is called a \emph{branching point}. The
branching points are isolated, hence there are finitely many, say
$n$, of them. Removing the branching points in $\Sigma$ and all
their pre-images in $\Sigmatil$, the restriction of $f$ gives a
genuine cover, whose degree we will denote by $d$. The collection of
local degrees at the preimages of the $i$-th branching point of
$\Sigma$ is a partition $\Pi_i$ of $d$, namely a set of positive
integers summing up
to $d$. Let us define $\ell(\Pi_i)$ as the length of this
partition, $\Pi=\{\Pi_1,\ldots,\Pi_n\}$ and
$\ell(\Pi)=\ell(\Pi_1)+\ldots+\ell(\Pi_n)$. The whole information on
the branched cover $f$ will be summarized by the symbol
$$f:\Sigmatil\myto^{d:1}_{\Pi}\Sigma.$$
The following is known:
\begin{prop}\label{nec:cond:prop}
If $f:\Sigmatil\myto^{d:1}_{\Pi_1,\ldots,\Pi_n}\Sigma$ is a surface branched cover then:
\begin{itemize}
\item[(1)]
$\chi(\Sigmatil)-\ell(\Pi)=d\cdot\big(\chi(\Sigma)-n\big)$;
\item[(2)]
$n\cdot d-\ell(\Pi)$ is even;
\item[(3)]
If $\Sigma$ is
orientable then $\Sigmatil$ is also orientable;
\item[(4)] If $\Sigma$ is non-orientable and $d$ is odd
then $\Sigmatil$ is also non-orientable;
\item[(5)] If
$\Sigma$ is non-orientable and $\Sigmatil$ is orientable then each
partition $\Pi_i$ of $d$ is a juxtaposition of two partitions of
$d/2$.
\end{itemize}
\end{prop}
Condition~(1) is the classical Riemann-Hurwitz formula.
Condition~(2) follows from (1)
in the orientable case, and it
is not too hard to establish in general, see~\cite{partI}.
Conditions~(3) and~(4) are obvious.
Condition~(5) is quite easy, see~\cite{partI} again; note that $d$ is even by
condition~(4).

\paragraph{Candidate surface branched covers}
Suppose we are given closed connected surfaces $\Sigmatil$ and
$\Sigma$, integers $n\geqslant 0$ and $d\geqslant 2$, and a collection
$\Pi=(\Pi_1,\ldots,\Pi_n)$ of $n$ partitions of $d$.
We will say that these data define a \emph{candidate surface branched cover}, and
we will associate to them the symbol
$$\Sigmatil\argdotsto{d:1}{\Pi_1,\ldots,\Pi_n}\Sigma$$
if the conditions of Proposition~\ref{nec:cond:prop} hold.

\begin{rem}
\emph{The symbol $\big(\Sigmatil,\Sigma,n,d,(d_{ij})\big)$ and the name
\emph{compatible branch datum} are used in \cite{partI,partII} instead
of the terminology of ``candidate covers'' we will use here.}
\end{rem}

A candidate surface branched cover $\Sigmatil\argdotsto{d:1}{\Pi}\Sigma$
will be called \emph{realizable} if there
is an actual cover $f:\Sigmatil\myto^{d:1}_{\Pi}\Sigma$
matching it, and \emph{exceptional} otherwise.
A classical problem dating back to Hurwitz~\cite{Hurwitz}
asks which candidate covers are realizable and which are exceptional.
Many authors contributed to it, as we will now outline, but the problem
in full generality is still unsolved.

\paragraph{Known results} A thorough description of the partial
solutions to the Hurwitz problem obtained over the time was given
in~\cite{partI}. Here we restrict ourselves to the main results. We
start with the following theorem of Husemoller~\cite{Husemoller},
proved also in~\cite{EKS}:

\begin{thm}
A candidate surface branched cover $\Sigmatil\argdotsto{d:1}{\Pi}\Sigma$
is realizable if $\Sigma$ is orientable and
$\chi(\Sigma)\leqslant 0$.
\end{thm}

We next have the following result which combines theorems of Ezell~\cite{Ezell}
and Edmonds-Kulkarny-Stong~\cite{EKS}:

\begin{thm}
A candidate surface branched cover $\Sigmatil\argdotsto{d:1}{\Pi}\Sigma$
is realizable if either $\Sigma$ is non-orientable and $\chi(\Sigma)\leqslant 0$,
or $\Sigma$ is the projective plane and $\Sigmatil$ is non-orientable.
\end{thm}

According to these results the problem of the realizability of
$\Sigmatil\argdotsto{d:1}{\Pi}\Sigma$
remains open only if $\Sigmatil$ is orientable and $\Sigma$ is either the sphere $S^2$
or the projective plane. However other results in~\cite{EKS}
allow to reduce the latter case to the former one. For this reason
we will restrict in the rest of this paper to the case $\Sigma=S^2$.
Many exceptions are known in this case, the easiest of which is
$$S^2\argdotstoter{4:1}{(2,2),(2,2),(3,1)}S^2,$$
but the general pattern of realizability and exceptionality remains elusive.
The following conjecture suggesting connections with number-theoretic facts
was however proposed in~\cite{EKS}:
\begin{conj}\label{prime:conj}
If $\Sigmatil\argdotsto{d:1}{\Pi}S^2$ is a
candidate surface branched cover and the
degree $d$ is a prime number then the candidate is realizable.
\end{conj}
Its validity was recently supported by the results and computer
experiments of Zheng~\cite{Zheng} and by the results
of~\cite{partII} and~\cite{PaPe}.

\bigskip

In the rest of this paper we will review five different techniques
employed to attack the Hurwitz problem, using them we will give five independent
proofs of the exceptionality of $S^2\argdotstoter{4:1}{(2,2),(2,2),(3,1)}S^2$,
and we will state the main results they have led to.
We address the reader to~\cite{partI,partII,PaPe} for more details.

\section{Permutations}\label{permu:section}
Hurwitz himself already showed that the problem of realizability of
a candidate surface branched cover can be reformulated, using the
notion of monodromy, in terms of permutations.
We will first provide an example of how this works and then we will
sketch the general technique and the main achievements obtained using it.

\paragraph{Sample application} We will now give our first proof of the
exceptionality of the candidate surface branched cover
$S^2\argdotstoter{4:1}{(2,2),(2,2),(3,1)}S^2$. Suppose, on the
contrary, that there exists a branched cover
$f:S^2\rightarrow S^2$ realizing it. Let $S^2_n$
denote the surface obtained from the sphere $S^2$ by removing $n$ open discs
with disjoint closures. Clearly, $f$ induces a genuine cover
$S^2_6\rightarrow S^2_3$ such that each boundary component of
$S^2_3$ is covered by two boundary components of $S^2_6$, and
the degrees of the restrictions to these components are (2,2),
(2,2), and (3,1).

Choose two simple proper disjoint arcs $\epsilon_1$,
$\epsilon_2$ cutting $S^2_3$ into a 2-disc $\Delta$, and denote by
$\epsilon_{i,\pm}$ the arcs in the boundary of $\Delta$ corresponding to
$\epsilon_i$. Now note that a genuine cover over a disc
is given by a disjoint union of discs, each of which is mapped homeomorphically to
the target. Therefore each degree-4 cover over
$S^2_3$ can be reconstructed in the following way. We first take
the disjoint union of $4$ copies $(\Delta^{(h)})_{h=1}^4$ of
$\Delta$, with the corresponding arcs $\epsilon^{(h)}_{i,\pm}$,
and we choose two permutations $\theta_1,\theta_2\in\permu_4$. Then
we glue each $\epsilon^{(h)}_{i,-}$ to
$\epsilon^{(\theta_i(h))}_{i,+}$ and we project to $S^2_3$ in the obvious way.

It is easy to show that the
degrees of the covers over the boundary components of $S^2_3$
are then given by the lengths of the cycles of $\theta_1$, $\theta_2$, and
$\theta_1\myprod\theta_2$. (Moreover the covering surface is
connected if and only if $\theta_1$ and $\theta_2$ generate a
subgroup of $\permu_4$ acting transitively on $\{1,\ldots,4\}$, but we will not
need this fact here).

We conclude that if the candidate
$S^2\argdotstoter{4:1}{(2,2),(2,2),(3,1)}S^2$ is realizable
then there exist permutations
$\theta_1,\theta_2\in\permu_4$ such that the lengths of the cycles of $\theta_1$,
$\theta_2$ and $\theta_1\myprod\theta_2$ are respectively $(2,2)$, $(2,2)$, and $(3,1)$.
However the set of all elements of $\permu_4$ with cycles of lengths $(2,2)$,
together with the trivial permutation, is a subgroup of $\permu_4$, whence
a contradiction and the conclusion that $S^2\argdotstoter{4:1}{(2,2),(2,2),(3,1)}S^2$
is exceptional.

\paragraph{General technique} The above argument can be generalized
to obtain the following result of Hurwitz~\cite{Hurwitz}. Again
$\Sigma_n$ denotes $\Sigma$ with $n$ discs removed. We confine
ourselves to the statement for the orientable case, a small
additional but rather technical condition being required in the
general case.

\begin{thm}\label{Hurwitz:method:thm}
A candidate surface branched cover
$\Sigmatil\argdotstobis{d:1}{\Pi_1,\ldots,\Pi_n}\Sigma$ with
orientable $\Sigma$ is realizable if and only if there exists a
representation $\theta:\pi_1(\Sigma_n)\rightarrow \permu_d$ such
that:
\begin{enumerate}
\item The image of $\theta$ acts transitively on $\{1,\ldots,d\}$;
\item The image under $\theta$ of the element of $\pi_1(\Sigma_n)$
corresponding to the $i$-th boundary component of $\Sigma_n$ has
cycles of lengths $\;\Pi_i$.
\end{enumerate}
\end{thm}

It is precisely this technique that was
employed in~\cite{Husemoller,EKS,Ezell} and led to the
results obtained therein and stated in the Introduction.
The proofs are not elementary at all, but the underlying philosophy is easy to
explain. The main restrictions imposed by
Theorem~\ref{Hurwitz:method:thm} are on the images under $\theta$ of the
peripheral elements of $\pi_1(\Sigma_n)$. But if
$\chi(\Sigma)\leqslant 0$ then $\pi_1(\Sigma_n)$ is \emph{not} generated by
these peripheral elements, so one has more flexibility for the choice of $\theta$.

\section{Dessins d'enfants}\label{dessins:section}
The notion of dessin d'enfant was introduced by Grothendieck~\cite{Groth}
for matters related to the Hurwitz problem but not strictly equivalent to it.
Again before giving the general definition we show in a concrete case
how the notion works.

\paragraph{Sample application} Our second proof of the exceptionality of the
sample candidate $S^2\argdotstoter{4:1}{(2,2),(2,2),(3,1)}S^2$ is based
on the following observation. Suppose that some $S^2\argdotstoter{4:1}{(2,2),(2,2),\Pi_3}S^2$
is realized by a map $f$. Denote the branching points by $x_1,x_2,x_3$, with
the partitions $(2,2)$ associated to $x_1$ and $x_2$. Let $\alpha$ be a
simple arc in the base $S^2$ that joins $x_1$ to $x_2$ and avoids
$x_3$. Then $f^{-1}(\alpha)$ is a graph $D$ in the covering
$S^2$ whose vertices are the pre-images of $x_1$ and $x_2$.
However the valence of any such vertex $v$ is the local degree of
$f$ at $v$, which is $2$, so $D$ is a union of circles.

\begin{figure}
\centering
\includegraphics[width=10cm]{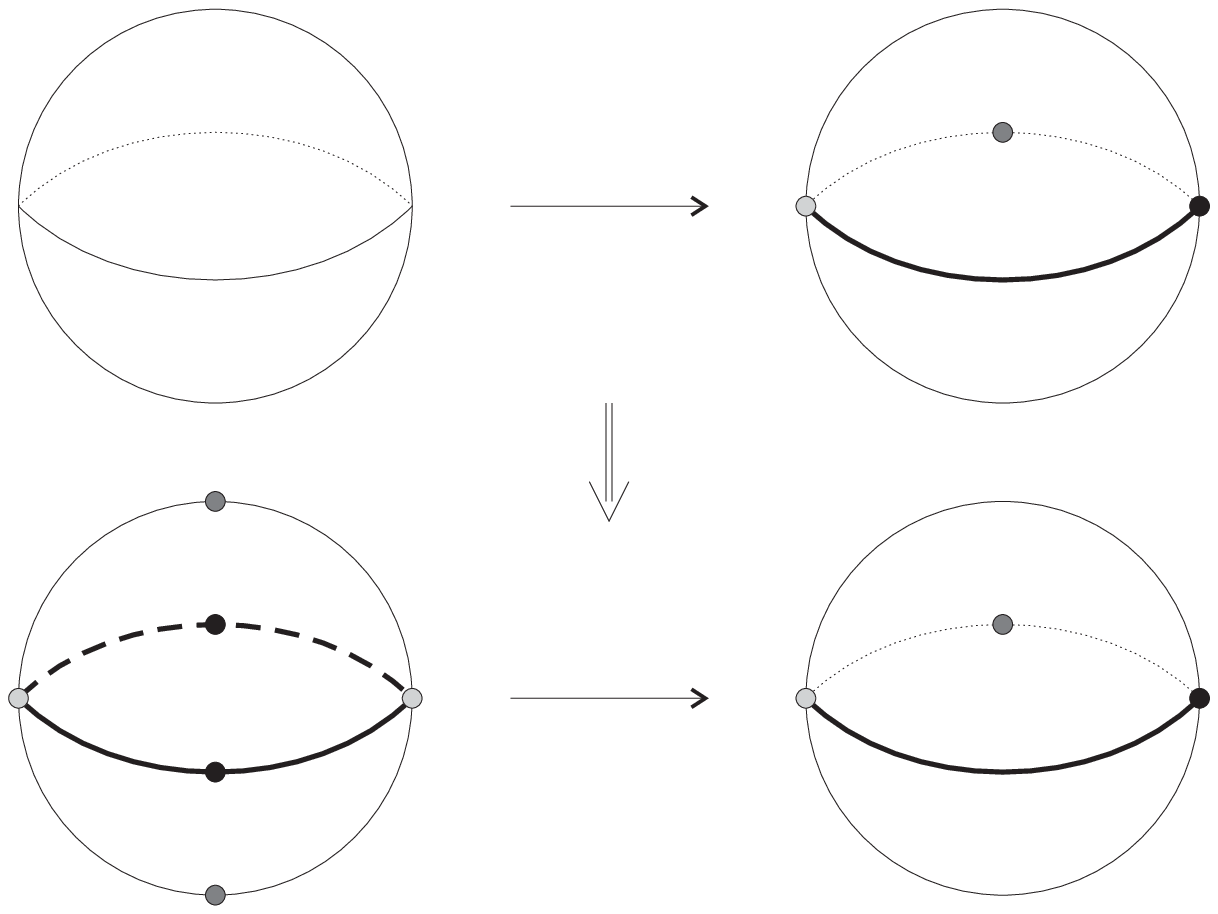}
\mycap{Dessin d'enfant for a cover of degree 4 with two partitions (2,2)\label{dessin:sample:fig}}
\end{figure}

Now note that the complement of $\alpha$ in the base $S^2$ is an open disc containing
only one branching point. This implies that the complement of $D$ in
the covering $S^2$ is a union of two discs,
each containing one element of $f^{-1}(x_3)$, so $D$ is a single circle (with 4 vertices,
even if one cannot see them).
Moreover the local degree at the element of $f^{-1}(x_3)$ contained in a disc
is half the number of vertices of $D$ that the disc is incident to, therefore it is $2$.
This shows that $\Pi_3$ is forced to be $(2,2)$, so
$S^2\argdotstoter{4:1}{(2,2),(2,2),(3,1)}S^2$
is exceptional. Our argument is pictorially illustrated in Figure~\ref{dessin:sample:fig}.

\paragraph{General technique}
Grothendieck's original dessins d'enfants~\cite{Groth,Wolfart} arise
precisely as in the argument we have just given, when one considers
a branched cover $f:\Sigmatil\to S^2$ with three branching points and one
defines the dessin $D$ as $f^{-1}(\alpha)$, where $\alpha$ is an arc joining
two of the branching points and avoiding the third one. This technique was
generalized in~\cite{partI} to the case of an arbitrary number of
branching points, as we will now explain.

\begin{defn}
\emph{A \emph{dessin d'enfant} on $\Sigmatil$ is a graph
$D\subset\Sigmatil$ where:
\begin{enumerate}
\item For some $n\geqslant 3$ the set of vertices of $D$ is split
as $V_1\sqcup\ldots\sqcup V_{n-1}$ and the set of edges of $D$ is
split as $E_1\sqcup\ldots\sqcup E_{n-2}$;
\item For $i=1,\ldots,n-2$
each edge in $E_i$ joins a vertex of $V_i$ to one of $V_{i+1}$;
\item For $i=2,\ldots,n-2$ any vertex of $V_i$ has even valence and
going around the vertex one alternatively encounters edges from $E_{i-1}$ and
from $E_{i}$;
\item $\Sigmatil\setminus D$ consists of open
discs.
\end{enumerate}}
\end{defn}

The following was established in~\cite{partI}.

\begin{prop}\label{dessins:covers:prop}
The realizations of a
candidate surface branched cover $\Sigmatil\argdotstobis{d:1}{\Pi_1,\ldots,\Pi_n}\Sigma$
correspond to the dessins
d'enfants $D\subset\Sigmatil$ with the set of vertices split as
$V_1\sqcup V_2\sqcup\ldots\sqcup V_{n-1}$ such that:
\begin{itemize}
\item For $i=1$ and
$i=n-1$ the vertices in $V_i$ have valences $\;\Pi_i$;
\item For $i=2,\ldots,n-2$ the vertices in
$V_i$ have valences $2\myprod \Pi_i$;
\item The numbers of vertices (with multiplicity) that
the discs in
$\Sigmatil\setminus D$ are incident to are $2(n-2)\myprod\Pi_n$.
\end{itemize}
\end{prop}

\paragraph{Some results}
Using Proposition~\ref{dessins:covers:prop} one can analyze the realizability
of several infinite series of candidate surface branched covers. For instance
the following results were establishes in~\cite{partI}:

\begin{prop}\label{non-realiz:53:prop}
Let $d\geqslant 8$ be even and consider a candidate surface branched cover of
the form
$$\Sigmatil\argdotstoqua{d:1}{(2,\ldots,2), (5,3,2,\ldots,2),\Pi_3}S^2.$$
\begin{itemize}
\item If $\Sigmatil$ is the torus, whence $\ell(\Pi_3)=2$, the candidate is
realizable if and only if $\;\Pi_3\neq(d/2,d/2)$;
\item If
$\Sigmatil$ is the sphere $S^2$, whence $\ell(\Pi_3)=4$, the candidate is realizable if and
only if $\;\Pi_3$ does not have the form
$(a,a,b,b)$ or $(3a,a,a,a)$ for $a,b\in\matN$.
\end{itemize}
\end{prop}

\begin{prop}\label{23:nonex:prop}
Let $d$ be even. Then a candidate surface branched cover of one of
the forms
$$S^2\argdotstoqua{d:1}{(2,\ldots,2),(3,3,2,\ldots,2),\Pi_3}S^2,\qquad
S^2\argdotstoqua{d:1}{(2,\ldots,2),(3,2,\ldots,2,1),\Pi_3}S^2.$$
is realizable if and only if the largest element of $\;\Pi_3$ is not $d/2$.
\end{prop}

The proofs of these propositions are not too difficult,
the idea being that a vertex of valence 2 is actually not a true vertex
of a graph. So in all cases one has to deal with a small number
of topological types of dessins d'enfant, and a slightly larger number
of inequivalent embeddings of these graphs in the relevant surface
$\Sigmatil$. Then one has to analyze how the (invisible) valence-2 vertices
are placed along the graph and to analyze what partitions $\Pi_3$
arise. A somewhat more complicated argument yields the following result.

\begin{thm}\label{excep:by:fixpoints:thm}
Suppose that a candidate surface branched cover
$S^2\argdotstobis{d:1}{\Pi_1,\ldots,\Pi_n}S^2$ is realizable and each entry
of $\;\Pi_1$ and $\;\Pi_2$ is a multiple of some
$k$ with $1<k<d$, so also $d$ is.
Then each entry of each $\;\Pi_i$ with $i\geqslant 3$
is at most $d/k$.
\end{thm}

\section{Checkerboard graphs}\label{checker:section}
This technique was introduced by Bar\'anski~\cite{Baranski} for the
case of candidate covers $S^2\dotsto S^2$ and extended to the case
$\Sigmatil\dotsto S^2$ in~\cite{partII}. As above we start with an
example and then we outline the general method.

\begin{figure}
\centering
\includegraphics[width=11cm]{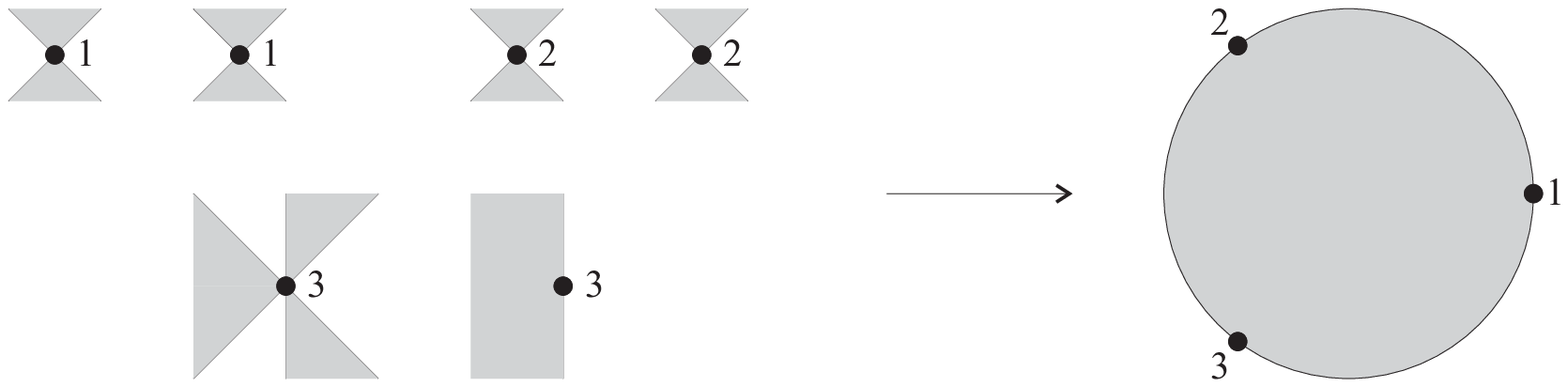}
\mycap{Vertices of a
checkerboard graph for a degree-4 cover of $S^2$
with partitions (2,2), (2,2), and (3,1)\label{checker:1:sample:fig}}
\end{figure}

\paragraph{Sample application} Here we present the third proof of
the exceptionality of $S^2\argdotstoter{4:1}{(2,2),(2,2),(3,1)}S^2$.
Suppose that the candidate is realized by some map $f$. Let us
identify $S^2$ with $\matC\cup\{\infty\}$ and assume that the
branching points are the third roots of $1$ in $\matC$, labelled
$1$, $2$ and $3$. Let $\Delta$ be the unit disc in $\matC$ and let
us colour it black and its complement white. Then
$f^{-1}(\partial\Delta)$ is a graph with 6 vertices, of valences
$4$, $4$, $4$, $4$, $6$, and $2$, each with one of the labels $1$,
$2$ or $3$ equal to that of its image. Moreover the complement of
the graph consists of 4 black and 4 white discs (because the degree
of $f$ is 4). Each disc is actually a (curvilinear) triangle, and
whenever two triangles share an edge they have different colours (as
in a checkerboard, whence the name of the technique). See
Figure~\ref{checker:1:sample:fig}.

\begin{figure}
\centering
\includegraphics[width=9cm]{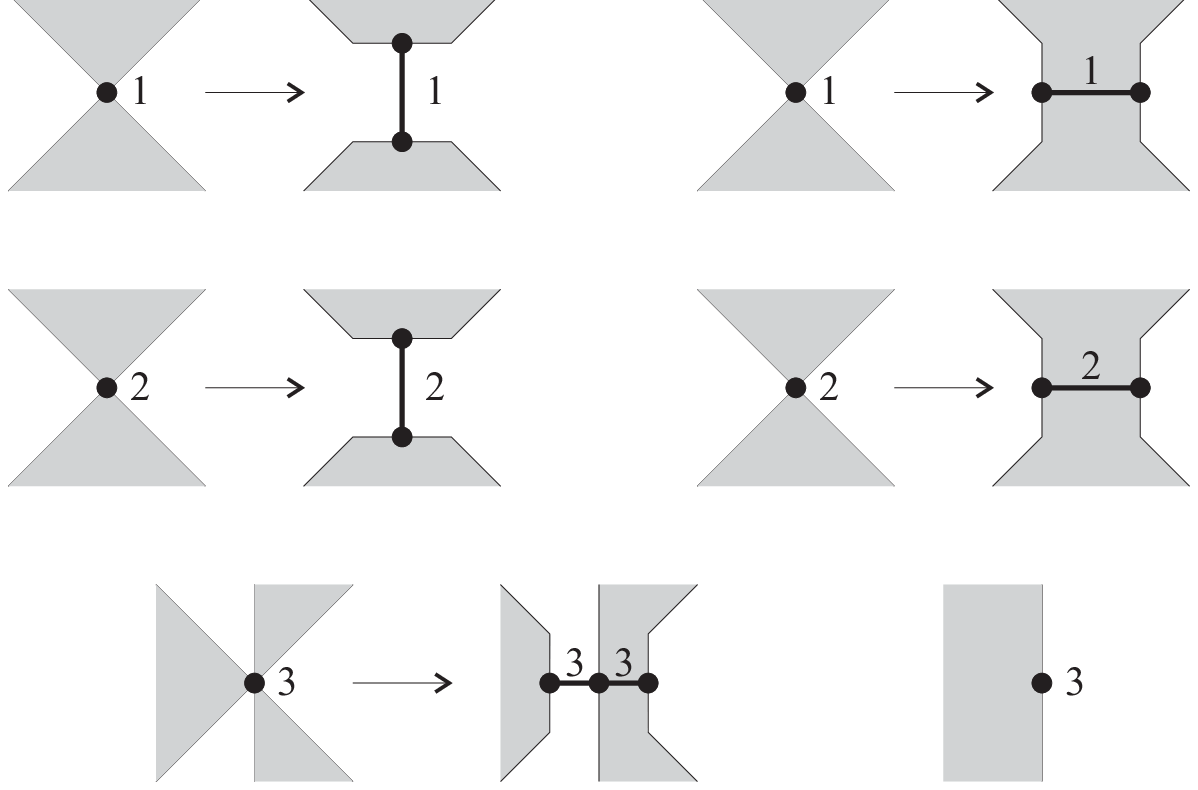}
\mycap{Moves employed to merge the triangles\label{checker:2:sample:fig}}
\end{figure}

We now modify our graph by merging together all the black
triangles into one black disc, and all the white triangles into
one white disc. This is done by moves such as those of
Figure~\ref{checker:2:sample:fig}.
The moves are applied successively and in such a way that each
time the pairs of local germs of discs merged together belong to
different global discs. In addition to performing this merging, as
also shown in Figure~\ref{checker:2:sample:fig}, we insert some
vertices and arcs labelled $1$, $2$ or $3$, with obvious
conventions. Note that the set of moves to apply is not unique
(for instance, a different set could
perhaps give rise to two arcs labelled $1$ both contained in the
black merged disc, or similar other variations), but the number of
actually distinct possibilities is finite and small.

\begin{figure}
\centering
\includegraphics[width=10cm]{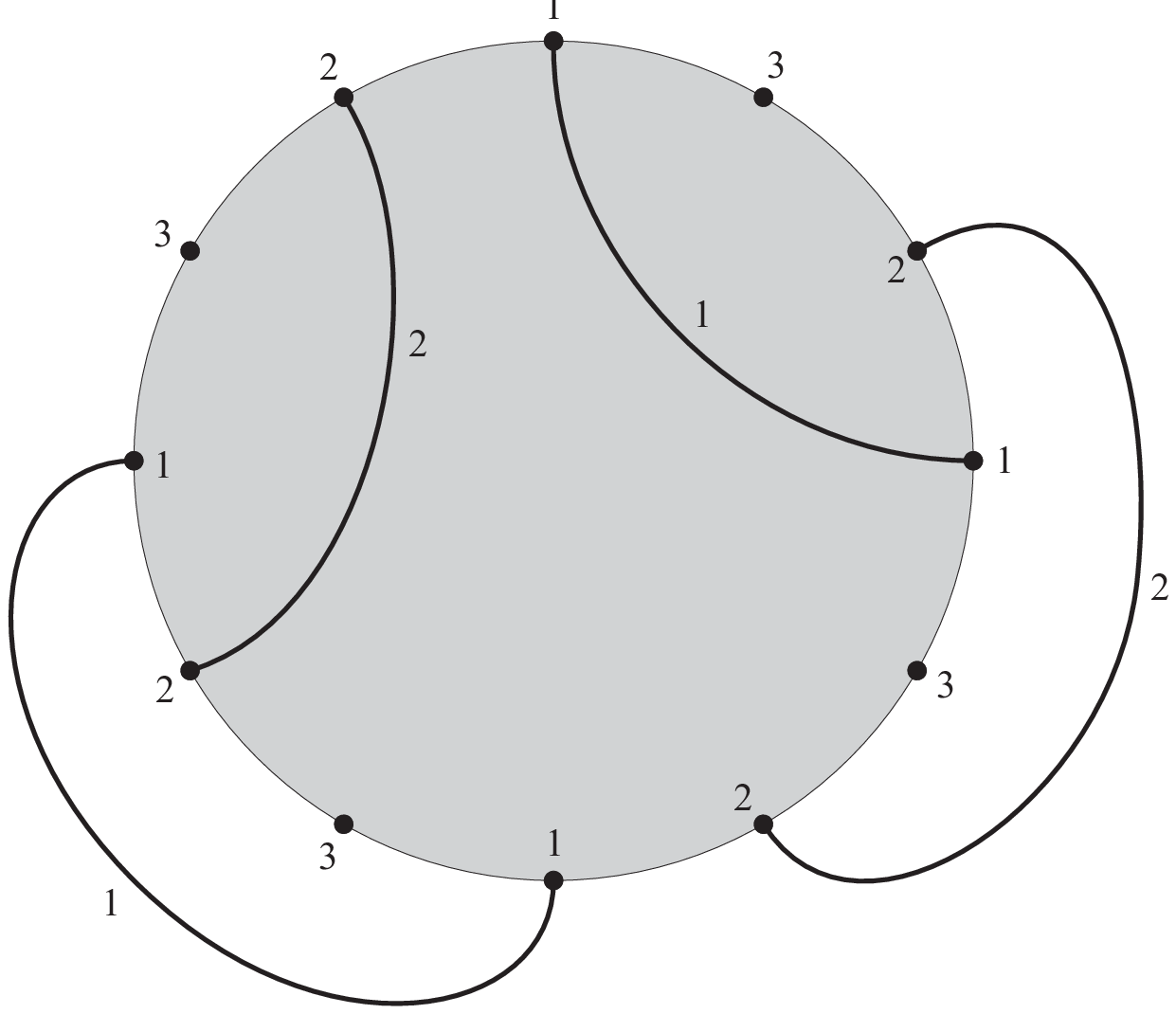}
\mycap{Merged black and white discs with labelled arcs\label{checker:3:sample:fig}}
\end{figure}

It is not difficult to prove that, after all the merging, on the
boundary of the single black disc there are 12 vertices, with labels
$1,2,3$ repeated cyclically 4 times. In addition there are two arcs
labelled $1$ with vertices $1$ at the ends, two similar arcs
labelled $2$, and an arc labelled $3$ with two vertices labelled $3$
at the ends and one in the middle. There are a few different
possibilities for the background colours of the arcs, and by examining
all of them one sees that this triple actually cannot exist.
One of the situations to consider is shown as an example in
Figure~\ref{checker:3:sample:fig},
where after insertion of the arcs labelled $1$ and $2$ it is
impossible to draw the arc labelled $3$, because the arcs must be
mutually disjoint. This shows that the initial assumption about the
existence of $f$ was absurd, which implies once again that
$S^2\argdotstoter{4:1}{(2,2),(2,2),(3,1)}S^2$ is exceptional.

\paragraph{General technique} To describe the general form of this approach,
we need a preliminary notion. We call \emph{checkerboard graph}
in a surface $\Sigmatil$ a
finite $1$-subcomplex $G$ of $\Sigmatil$ whose
complement consists of open discs each bearing a color black or
white, so that each edge separates black from white.

The strategy used in the example extends as follows to the general
situation. Suppose we have a map realizing a candidate surface
branched cover $\Sigmatil\argdotstobis{d:1}{\Pi_1,\ldots,\Pi_n}S^2$.
Then we arrange the branching points to be the $n$-th roots of $1$
and we consider the graph $f^{-1}(\partial\Delta)$. Giving black
colour to the complementary discs mapped to $\Delta$, and white to
the other discs, we see that $f^{-1}(\partial\Delta)$ is a
checkerboard graph in $\Sigmatil$. Putting labels on the vertices
and performing moves as those in Figure~\ref{checker:2:sample:fig}
(each consisting in the merging of two discs and the insertion of a
labelled arc), we end up with a checkerboard graph whose complement
consists of a single black and a single white disc, together with a
collection of labelled vertices and disjoint trees satisfying a long
list of conditions in terms of the partitions $\Pi_1,\ldots,\Pi_n$.
The precise list would be too long to give here, so we address the
reader to~\cite{partII}. The main point is however that the whole
process is reversible, namely from a checkerboard graph and a family
of labelled vertices and trees satisfying the conditions
corresponding to some candidate branched surface cover
$\Sigmatil\argdotstobis{d:1}{\Pi_1,\ldots,\Pi_n}S^2$ one can
construct a map realizing the candidate. Thus checkerboard graphs
give a necessary and sufficient criterion for realizability.

\paragraph{Results}
The following main theorem was established in~\cite{partII}
by means of the checkerboard graphs realizability criterion.
The proof was carried out using some inductive constructions.

\begin{thm}
Consider a candidate surface branched cover of the form
$$\Sigmatil\argdotstoter{d:1}{(d-2,2),\Pi_2,\Pi_3}S^2.$$
\begin{itemize}
\item If $\Sigmatil$ is the sphere $S^2$
then the candidate is exceptional if it has one of the forms
$$
S^2\argdotstoqui{2k:1}{(2k-2,2),(2,\ldots,2),(k+1,1,\ldots,1)}S^2
,\qquad
S^2\argdotstoqua{2k+2:1}{(2k,2),(2,\ldots,2),(2,\ldots,2)}S^2
$$
for $k\geqslant 2$, and realizable otherwise. In particular, it is
always realizable if the degree $d$ is odd.
\item If $\Sigmatil$ is the torus $S^1\times S^1$ then the candidate is always realizable,
with the single exception of
$$S^1\times S^1\argdotstoter{6:1}{(4,2),(3,3),(3,3)}S^2.$$
\item If $\Sigmatil$ has genus at least $2$ then the candidate is always realizable.
\end{itemize}
\end{thm}

\section{Factorization}\label{factor:section}
A composition of branched covers is a branched cover, so if a
candidate can be expressed as a ``candidate composition'' of two
realizable covers then it is realizable. On the other hand if one
has a candidate cover and one can show that a map realizing it, if
any, should be the composition of two maps one of which realizes an
exceptional cover, one deduces that the original cover is
exceptional too.

\begin{figure}
\centering
\includegraphics[width=10cm]{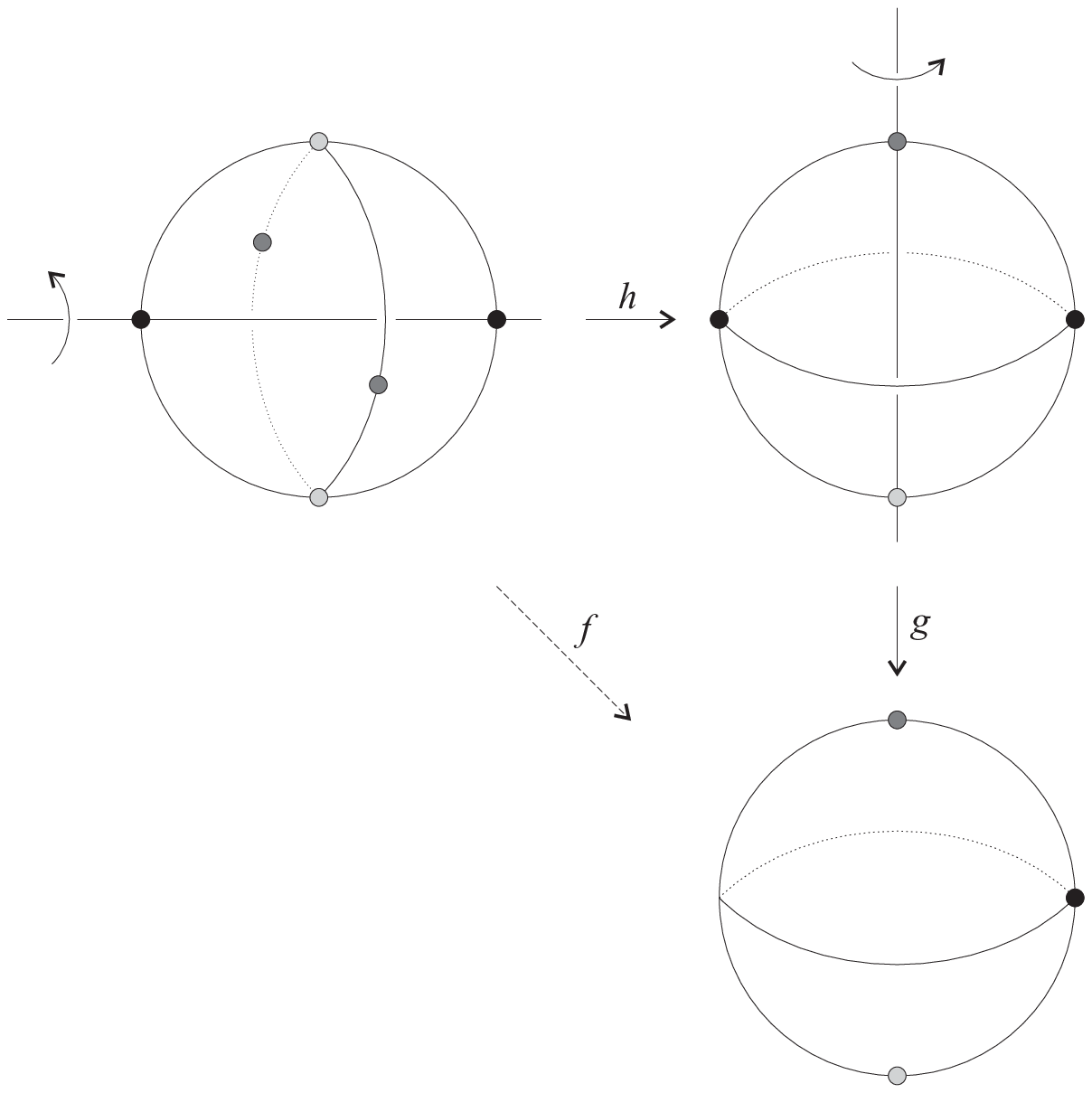}
\mycap{Stepwise construction of a degree-4 cover\label{factor:sample:fig}}
\end{figure}

\paragraph{Sample application}
We will now give our fourth proof of the exceptionality of
$S^2\argdotstoter{4:1}{(2,2),(2,2),(3,1)}S^2$.  The proof is actually
slightly incomplete, but see below. Suppose we want to construct
a map $f$ realizing some candidate cover
$S^2\argdotstobis{4:1}{(2,2),(2,2),\Pi_3}S^2$ and we decide to
do this stepwise. Let us denote by $x_1,x_2,x_3$ the
would-be branching points of $f$.
We start by realizing the local degrees
$2$ at $x_1$ and $x_2$ via a degree-2 cover $g:S^2\to S^2$,
but there is essentially just one of them, whose only non-trivial automorphism
is a rotation of angle $\pi$ through two antipodal points, see
the right portion of Figure~\ref{factor:sample:fig}.
Note that $g^{-1}(x_3)$ consists of two points.
Now we want to find another degree-2 cover $h:S^2\to S^2$ such that $f=g\compo h$.
But we know that $\Pi_3$ must have length $2$, so the two points in $g^{-1}(x_3)$ must
be the branching points of $h$ (which coincides with $g$ up to automorphisms of $S^2$).
So the situation is that Figure~\ref{factor:sample:fig}, which implies that
$\Pi_3=(2,2)$.

This shows that no map $f$ realizing
$S^2\argdotstoter{4:1}{(2,2),(2,2),(3,1)}S^2$ can be constructed
using the stepwise method we have outlined. As we will see below,
one can show that the candidate satisfies properties ensuring that
if it is realizable then there exists a stepwise realization, so
again we conclude that it is exceptional.

\paragraph{General technique and results}
The missing portion of the previous argument is deduced from the following
result established in~\cite{partI}, together with its main consequence
stated soon after.

\begin{prop}\label{filtration:prop}
Let a candidate surface branched cover
$S^2\argdotstobis{d:1}{\Pi_1,\ldots,\Pi_n}S^2$
be realized by a map $f$ and suppose that
all the entries of $\;\Pi_1$ and $\;\Pi_2$ are even,
so $d$ also is. Then $f$ can be expressed as $f=g\compo h$ where $g$ is the
realizable cover $S^2\myto^{2:1}_{(2),(2)}S^2$
and $h$ realizes a candidate of the form
$$S^2\argdotstosex{d/2:1}{\frac 12\Pi_1,\frac 12\Pi_2,\Pi'_3,\Pi''_3,\ldots,\Pi'_n,\Pi''_n}S^2$$
where $\;\Pi_i$ is the juxtaposition of $\;\Pi'_i$ and $\;\Pi''_i$ for $i\geqslant 3$.
\end{prop}

\begin{thm}\label{even-deg_exceptions:thm}
If $d$ and all the entries of $\;\Pi_1$ and $\;\Pi_2$ are even and
the candidate surface branched cover
$$S^2\argdotstobis{d:1}{\Pi_1,\ldots,\Pi_n}S^2$$
is realizable then $\;\Pi_i$ is the juxtaposition of two partitions of $d/2$
for $i\geqslant 3$.
\end{thm}

On the realizability side, the following easy but useful fact
was also established in~\cite{partI} exploiting the idea of factorization.

\begin{thm}\label{divisible:exist:thm}
Consider a candidate surface branched cover
$$\Sigmatil\argdotstobis{d:1}{\Pi_1,\Pi_2,\Pi_3}S^2$$
and suppose there exists an odd number $p\geqslant 3$
dividing each entry in each $\Pi_i$. Then the candidate is realizable.
\end{thm}

\section{Geometric 2-orbifolds}\label{orb:section}
In this last section, before giving yet another proof of the exceptionality of
$S^2\argdotstoter{4:1}{(2,2),(2,2),(3,1)}S^2$, we will have to review some
general notions.

\paragraph{Orbifolds and orbifold covers}
A 2-orbifold $X=\Sigma(p_1,\ldots,p_n)$ is
a closed orientable surface $\Sigma$
with $n$ cone points of orders $p_i\geqslant 2$, at which $X$ has a singular
differentiable structure given by the quotient $\matC/_{\langle{\rm rot}(2\pi/p_i)\rangle}$.
Thurston~\cite{thurston:notes} defined an orbifold cover $f:\Xtil\myto^{d:1}X$
of degree $d$ as a map such that generic points have $d$ preimages, and $f$ is locally modelled on
functions of the form
$$\matC/_{\langle{\rm rot}(2\pi/\ptil)\rangle}\myto^{k:1}\matC/_{\langle{\rm rot}(2\pi/p)\rangle}$$
induced by the identity of $\matC$, where $p=k\cdot\ptil$. He also introduced the notion of orbifold
Euler characteristic
$$\chiorb\big(\Sigma(p_1,\ldots,p_n)\big)=\chi(\Sigma)-\sum_{i=1}^n\left(1-\frac1{p_i}\right),$$
designed so that if $f:\Xtil\myto^{d:1} X$ is an orbifold cover then
$\chiorb(\Xtil)=d\cdot\chiorb(X)$. He then proved that $2$-orbifolds are almost
always \emph{geometric}, that is:
\begin{itemize}
\item If $\chiorb(X)>0$ then $X$ is either \emph{bad} (not orbifold-covered by a surface) or
\emph{spherical}, namely the quotient of the metric 2-sphere $\matS^2$ under a finite isometric action;
\item If $\chiorb(X)=0$ (respectively, $\chiorb(X)<0$) then $X$ is
\emph{Euclidean} (respectively, \emph{hyperbolic}), namely
the quotient of the Euclidean plane $\matE^2$ (respectively, the hyperbolic plane $\matH^2$)
under a discrete isometric action.
\end{itemize}
Finally, he showed that any orbifold has an \emph{orbifold universal cover},
which easily implies the following:
\begin{lemma}\label{no:bad:on:good:lem}
If $\Xtil$ is bad and $X$ is good then there cannot exist any orbifold cover $\Xtil\to X$.
\end{lemma}

\paragraph{Induced orbifold covers}
The local model described above for an orbifold cover $f:\Xtil\to X$
can be viewed as the map $\matC\ni z\mapsto z^k\in\matC$, therefore
$f$ is also a branched cover $\Sigmatil\to\Sigma$ between the
surfaces underlying $\Xtil$ and $X$. However different $\Xtil\to X$
can give the same $\Sigmatil\to\Sigma$, because in the local model
one can arbitrarily replace $p$ by some $h\myprod p$ and $\ptil$ by
$h\myprod \ptil$. On the other hand, any $\Sigmatil\to\Sigma$ has an
associated ``easiest'' $\Xtil\to X$, where the cone orders are
chosen as small as possible. This carries over to \emph{candidate}
covers, as we will now spell out. Consider a candidate surface
branched cover
$$\Sigmatil\argdotstosex{d:1}{(d_{11},\ldots,d_{1m_1}),\ldots,(d_{n1},\ldots,d_{nm_n})}\Sigma$$
and define
$$\begin{array}{ll}
p_i={\rm l.c.m.}\{d_{ij}:\ j=1,\ldots,m_i\},\quad & p_{ij}=p_i/d_{ij},\phantom{\Big|} \\
X=\Sigma(p_1,\ldots,p_n), & \Xtil=\Sigmatil\big((p_{ij})_{i=1,\ldots,n}^{j=1,\ldots,m_i}\big)
\end{array}$$
where ``l.c.m.'' stands for ``least common multiple.''
Then we have an induced candidate $2$-orbifold cover $\Xtil\dotsto^{d:1}X$
satisfying the condition $\chiorb(\Xtil)=d\myprod\chiorb(X)$, which is easily
deduced from the Riemann-Hurwitz condition.

\paragraph{Sample application}
If the candidate $S^2\argdotstoter{4:1}{(2,2),(2,2),(3,1)}S^2$ were realizable, then the induced
candidate orbifold cover
$$S^2(3)\dotsto S^2(2,2,3)$$
would also
be realizable, which is impossible by Lemma~\ref{no:bad:on:good:lem}
because $S^2(3)$ is bad and $S^2(2,2,3)$ is good.

\paragraph{General technique}
If a candidate orbifold cover $\Xtil\dotsto^{d:1}X$ is complemented with the
instructions of which cone points of $\Xtil$ should be mapped to which cone points
of $X$, one can reconstruct a unique candidate surface branched
cover $\Sigmatil\argdotsto{d:1}{\Pi}\Sigma$, so one can fully switch to the
viewpoint of candidate orbifold covers. The $X$'s such
that $\chiorb(X)\geqslant 0$, namely the bad, spherical or Euclidean $X$'s, are
easily listed, and in addition on most of them the geometric structure (if any)
is \emph{rigid}. With these facts in mind the following program appears to be natural:
\begin{itemize}
\item Determine all the candidate surface branched covers inducing candidate orbifold covers
$\Xtil\dotsto X$ with $\chiorb(X)\geqslant 0$;
\item For any such a candidate, analyze its realizability using geometric methods, namely
Lemma~\ref{no:bad:on:good:lem} when $\Xtil$ is bad, or the fact that a map realizing the cover
can be lifted to an isometry $\matS^2\to\matS^2$ (or $\matE^2\to\matE^2$) when $\Xtil$ and $X$ are
spherical (or Euclidean; it turns out that $X$ is never bad).
\end{itemize}

The program has been fully carried out in~\cite{PaPe}, leading to the following main results:

\begin{thm}\label{pos:chi:main:thm}
Let a candidate surface branched cover $\Sigmatil\dotsto^{d:1}_\Pi\Sigma$
induce a candidate $2$-orbifold cover $\Xtil\dotsto^{d:1}X$ with $\chiorb(X)>0$.
Then $\Sigmatil\dotsto^{d:1}_\Pi\Sigma$ is exceptional if and only if
$\Xtil$ is bad and $X$ is spherical. All exceptions occur with
non-prime degree.
\end{thm}

\begin{thm}\label{244:main:thm}
Suppose $d=4k+1$ for $k\in\matN$. Then
$$S\argdotstosex{d:1}
{(\underbrace{{\scriptscriptstyle 2,\ldots,2}}_{2k},1),
(\underbrace{{\scriptscriptstyle 4,\ldots,4}}_{k},1),
(\underbrace{{\scriptscriptstyle 4,\ldots,4}}_{k},1)}S$$ is a
candidate surface branched cover, and it is realizable if and only
if $d$ can be expressed as $x^2+y^2$ for some $x,y\in\matN$.
\end{thm}

\begin{thm}\label{236:main:thm}
Suppose $d=6k+1$ for $k\in\matN$. Then
$$S\argdotstosex{d:1}
{(\underbrace{{\scriptscriptstyle 2,\ldots,2}}_{3k},1),
(\underbrace{{\scriptscriptstyle 3,\ldots,3}}_{2k},1),
(\underbrace{{\scriptscriptstyle 6,\ldots,6}}_{k},1)}S$$
is a candidate surface branched cover and
it is realizable if and only if $d$ can be expressed as $x^2+xy+y^2$ for some $x,y\in\matN$.
\end{thm}

We conclude by noting that the last two results provide strong supporting evidence
for Conjecture~\ref{prime:conj}, because of the following facts:
\begin{itemize}
\item A prime number of the form $4k+1$ can always be expressed as
$x^2+y^2$ for $x,y\in\matN$ (Fermat);
\item A prime number of the form $6k+1$
can always be expressed as
$x^2+xy+y^2$ for $x,y\in\matN$ (Gauss);
\item The integers that can be expressed as $x^2+y^2$ or as
$x^2+xy+y^2$ with $x,y\in\matN$ have asymptotically zero density in $\matN$.
\end{itemize}
This means that a candidate cover in any of these two statements is
``exceptional with probability 1,'' even though it is realizable
when its degree is prime. Note also that it was shown in~\cite{EKS}
that establishing Conjecture~\ref{prime:conj} in the special case of
three branching points would imply the general case.

\vspace{1cm}

\noindent Dipartimento di Matematica Applicata\\
Via Filippo Buonarroti, 1C\\
56127 PISA -- Italy\\
{\tt pervova@csu.ru}\\
{\tt petronio@dm.unipi.it}

\end{document}